\documentclass{llncs}

\usepackage{amsmath,amssymb,graphicx,calc,bbold}

\title{Low-Dimensional Faces of Random 0/1-Polytopes%
  \thanks{This work was done while the author was a member
    of the \emph{Mathematical Sciences Research Institute} at Berkeley, CA, 
    during Oct/Nov 2003.
  }
} 
\titlerunning{Low-Dimensional Faces of Random 0/1-Polytopes} 
\toctitle{Low-Dimensional Faces of Random 0/1-Polytopes} 
\author{Volker Kaibel}
\authorrunning{Volker Kaibel} 
\tocauthor{Volker Kaibel (TU Berlin)}
\institute{%
  DFG Research Center \emph{Mathematics for key technologies}\\
  TU Berlin MA 6--2\\
  Stra\ss e des 17. Juni 136\\
  10623 Berlin\\ 
  Germany\\
  \email{kaibel@math.tu-berlin.de}
}

\DeclareMathOperator{\N}{\mathbb{N}}
\DeclareMathOperator{\R}{\mathbb{R}}
\DeclareMathOperator{\probOp}{\mathbb{P}\,}
\DeclareMathOperator{\expectOp}{\mathbb{E}\,}
\DeclareMathOperator{\convOp}{conv\,}
\DeclareMathOperator{\affOp}{aff\,}
\DeclareMathOperator{\cubeOp}{Q}
\DeclareMathOperator{\cubeVertsOp}{V}
\DeclareMathOperator{\interiorOp}{int\,}
\DeclareMathOperator{\littleOOp}{o}
\DeclareMathOperator{\bigOOp}{O}

\DeclareMathOperator{\boundaryOp}{\partial}
\newcommand{\prob}[1]{\probOp[{#1}]}
\newcommand{\probBig}[1]{\probOp\big[{#1}\big]}
\newcommand{\expect}[1]{\expectOp[{#1}]}
\newcommand{\boundary}[1]{\boundaryOp({#1})}
\newcommand{\conv}[1]{\convOp({#1})}
\newcommand{\cube}[1]{\cubeOp_{#1}}
\newcommand{\probCond}[2]{\probOp[{#1}\,|\,{#2}]}
\newcommand{\probCondBig}[2]{\probOp\big[{#1}\,\big|\,{#2}\big]}
\newcommand{\expectCondBig}[2]{\expectOp\big[{#1}\,\big|\,{#2}\big]}
\newcommand{\cubeVerts}[1]{\cubeVertsOp_{#1}}
\newcommand{\setdef}[2]{\{{#1}\,:\,{#2}\}}
\newcommand{\setdefBig}[2]{\big\{{#1}\,:\,{#2}\big\}}
\newcommand{\zeroVec}{\mathbb{0}}
\newcommand{\oneVec}{\mathbb{1}}
\newcommand{\littleO}[1]{\littleOOp({#1})}
\newcommand{\bigO}[1]{\bigOOp({#1})}

\newcommand{\iuar}{independently uniformly at random}
\newcommand{\uar}{uniformly at random}

\spnewtheorem{dfm}[theorem]{DFM}{\bfseries}{\itshape}

\begin{document}

\maketitle

% --------------------------------------------------
% abstract
% --------------------------------------------------

\begin{abstract}
  Let~$P$ be a random 0/1-polytope in~$\R^d$ with~$n(d)$ vertices, and
  denote by~$\varphi_k(P)$ the \emph{$k$-face density} of~$P$, i.e.,
  the quotient of the number of $k$-dimensional faces of~$P$ and
  $\binom{n(d)}{k+1}$. For each $k\ge 2$, we establish the existence
  of a sharp threshold for the $k$-face density and determine the
  values of the threshold numbers~$\tau_k$ such that, for all
  $\varepsilon>0$,
  $$
  \expect{\varphi_k(P)}\ =\ 
  \begin{cases}
    1-\littleO{1} & \text{ if $n(d)\le 2^{(\tau_k-\varepsilon)d}$ for all~$d$} \\   
    \littleO{1}   & \text{ if $n(d)\ge 2^{(\tau_k+\varepsilon)d}$ for all~$d$} \\   
  \end{cases}
  $$
  holds for the expected value of~$\varphi_k(P)$. The threshold for
  $k=1$ has recently been determined in~\cite{KR03}.
  
  In particular, these results indicate that the high face densities
  often encountered in polyhedral combinatorics (e.g., for the
  cut-polytopes of complete graphs) should be considered more as a
  phenomenon of the general geometry of 0/1-polytopes than as a
  feature of the special combinatorics of the underlying problems.
\end{abstract}

% --------------------------------------------------
% Introduction
% --------------------------------------------------

\section{Introduction and Results}
\label{sec:intro}

Over the last decades, investigations of various special classes of
0/1-polytopes (convex hulls of sets of 0/1-points) have not only lead
to beautiful structural results on combinatorial optimization
problems, but also to powerful algorithms. Consequently, there has
been some effort to learn more about the general class of
0/1-polytopes (see~\cite{Zie00}). 

In the 1980's, e.g., several results on the graphs of 0/1-polytopes
have been obtained, most notably Naddef's proof~\cite{Nad89} showing
that they satisfy the Hirsch-conjecture. A quite spectacular
achievement in~2000 was B\'ar\'any and P\'or's theorem~\cite{BP01}
stating that random 0/1-polytopes (within a certain range of vertex
numbers) have super-exponentially (in the dimension) many facets.
Their proof is based on the methods developed in the early 1990's by
Dyer, F\"uredi, and McDiarmid~\cite{DFM92}, in order to show that the
expected volume of a random $d$-dimensional 0/1-polytope with~$n$
vertices drops from (almost) zero to (almost) one very quickly
with~$n$ passing the threshold $2^{(1-(\log e)/2)d}$.

While B\'ar\'any and P\'or's result sheds some light on the
highest-dimensional faces of random 0/1-polytopes, we investigate
their lower dimensional faces in this paper.  For a polytope~$P$
with~$n$ vertices and some $k\in[\dim P]$ (with
$[a]:=\{1,2,\dots,\lfloor a\rfloor\}$), we call
$$
\varphi_k(P)\ :=\ \frac{f_k(P)}{\binom{n}{k+1}}
$$
the \emph{$k$-face density} of~$P$, where $f_k(P)$ is the number of
$k$-dimensional faces of~$P$. Clearly, we have $0<\varphi_k(P)\le 1$,
and $\varphi_k(P)=1$ holds if and only if~$P$ is $(k+1)$-neighbourly
in the usual polytope theoretical sense (see, e.g., \cite{Zie95}).

The $1$-face density $\varphi_1(P)$ is the density of the graph
of~$P$. In this case, a threshold result for random 0/1-polytopes has
recently been obtained in~\cite{KR03}. However, for specific classes
of 0/1-polytopes, high $k$-face densities have been observed also for
larger values of~$k$. For example, the cut-polytopes of complete
graphs have $2$-face density equal to one (and thus, also $1$-face
density equal to one), i.e., every triple of vertices makes a
triangle-face (see \cite{BM86,DL97}). Note that the cut-polytopes of
complete graphs have $2^{\Theta(\sqrt{d})}$ vertices.

Here, we obtain that there is a sharp threshold for the $k$-face
density of random 0/1-polytopes for all (fixed)~$k$. The threshold
values nicely extend the results for $k=1$, while the proof becomes
more involved and needs a heavier machinery (the one developed in the
above mentioned paper by Dyer, F\"uredi, and McDiarmid). As a
pay-back, the proof, however, reveals several interesting insights
into the geometry of (random) 0/1-polytopes.

\subsection{Results}

Let us fix some $k\in\{1,2,\dots\}$, set $r:=k+1$, and
let $n\,:\,\N\rightarrow\N$ be a function (with $n(d)\in[2^d]$).

Define 
$$
\cubeVerts{d}\ :=\ \{0,1\}^d
\qquad\text{and}\qquad
\cube{d}\ :=\ [0,1]^d\ =\ \convOp\cubeVerts{d}\ ,
$$
and  consider the following two models of random 0/1-polytopes.

For the first one, choose~$W$ uniformly at random from
the $n(d)$-element subsets of~$\cubeVerts{d}$, and define $P_1:=\convOp
W$. This is the model referred to in the abstract.

For the second one, choose 
$S_1,\dots,S_r,X_1,\dots,X_{n(d)-r}\in\cubeVerts{d}$ \iuar, and define
$$
S:=\{S_1,\dots,S_r\}\ , \quad
X:=\{X_1,\dots,X_{n(d)-r}\}\ , \quad
P_2:=\conv{X\cup S}\ .
$$

The main part of the paper will be concerned with the proof of a
threshold result (Theorem~\ref{thm:main}) within the second model.
If, for some $\varepsilon>0$, $n(d)\le 2^{(\frac{1}{2}-\varepsilon)d}$
holds for all~$d$, then $S_1,\dots,S_r,X_1,\dots,X_{n(d)-r}$ are
pairwise different with high probability:
\begin{equation}
  \label{eq:intro:neu:1}
  \probBig{|S\cup X| = n(d)}\ =\ 1-\littleO{1}  
\end{equation}
This will allow us to deduce from Theorem~\ref{thm:main} the
threshold result within the first model promised in the abstract.

Throughout the paper, $\log(\cdot)$ and $\ln(\cdot)$ will denote the
binary and the natural logarithm, respectively.  For $0<\xi<1$, define
$$
h(\xi)\ :=\ \xi\log\frac{1}{\xi} +
            (1-\xi)\log\frac{1}{1-\xi}
$$
(i.e., $h(\cdot)$ is the binary entropy function).
Let us define
$$
H_r\ :=\ \frac{1}{2^r-2} 
         \sum_{i\in[r-1]}\binom{r}{i}h\big(\frac{i}{r}\big)
$$
and
$$
\Tilde{\tau}_r\ =\ 1-(1-2^{1-r})H_r\ .
$$
Note that we have $H_2=1$ and $0<H_r<1$ for $r\ge 3$.

\begin{theorem}
\label{thm:main}
  Let $r\in\{3,4,\dots\}$ and $\varepsilon>0$.
  \begin{enumerate}
  \item If $n(d)\le 2^{(\Tilde{\tau}_r-\varepsilon)d}$ holds for all~$d$, then we have
      $$
      \prob{\convOp S\text{ is a face of }P_2}\ =\ 1-\littleO{1}\ .
      $$
  \item If $n(d)\ge 2^{(\Tilde{\tau}_r+\varepsilon)d}$ holds for all~$d$, then we have
      $$
      \prob{P_2\cap \affOp S\text{ is a face of }P_2}\ =\ \littleO{1}\ .
      $$
  \end{enumerate}
\end{theorem}

From the evolution result on the density of the graphs of random
0/1-polytopes obtained in~\cite{KR03} one readily derives that the
statement of Theorem~\ref{thm:main} is also true for $r=2$ (note
$\Tilde{\tau}_2=\tfrac{1}{2}$).

Using Theorem~\ref{thm:main} (for $r\in\{2,3,\dots\}$), we can now prove the
main result of the paper, where for $k\in\{1,2,\dots\}$ we denote
$$
\tau_k\ :=\ \Tilde{\tau}_{k+1}\ =\ 1-(1-2^{-k})H_{k+1}\ .
$$

\begin{theorem}
  \label{thm:truemain}
  Let $k\in\{1,2,\dots\}$, $\varepsilon>0$, and $n:\N\rightarrow\N$ be
  any function. For each $d\in\N$, choose an $n(d)$-element subset~$W$
  of $\{0,1\}^d$ uniformly at random, and set $P:=\convOp W$. Then
  $$
  \expect{\varphi_k(P)}\ =\ 
  \begin{cases}
    1-\littleO{1} & \text{ if $n(d)\le 2^{(\tau_k-\varepsilon)d}$ for all~$d$} \\   
    \littleO{1}   & \text{ if $n(d)\ge 2^{(\tau_k+\varepsilon)d}$ for all~$d$} \\   
  \end{cases}
  $$
  holds for the expected $k$-face density of~$P$.
\end{theorem}

\begin{figure}[ht]
  \centering
  \includegraphics[height=4cm]{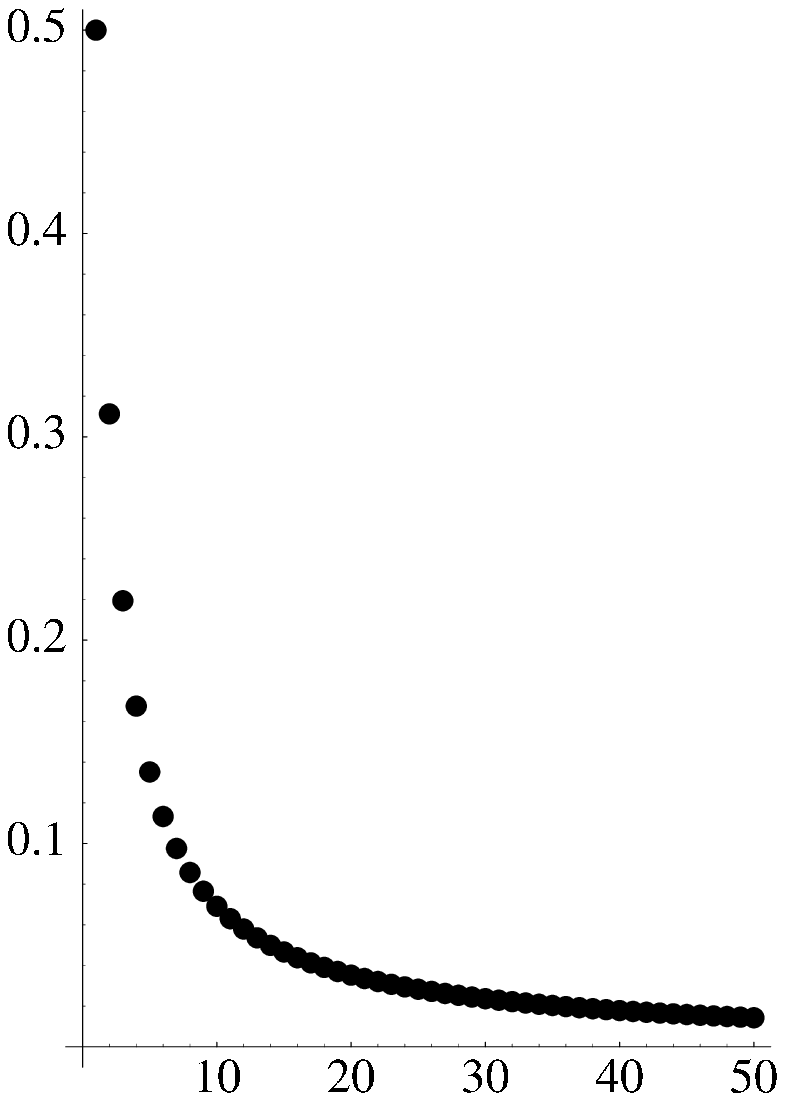}
  \hspace{2cm}
  \includegraphics[height=4cm]{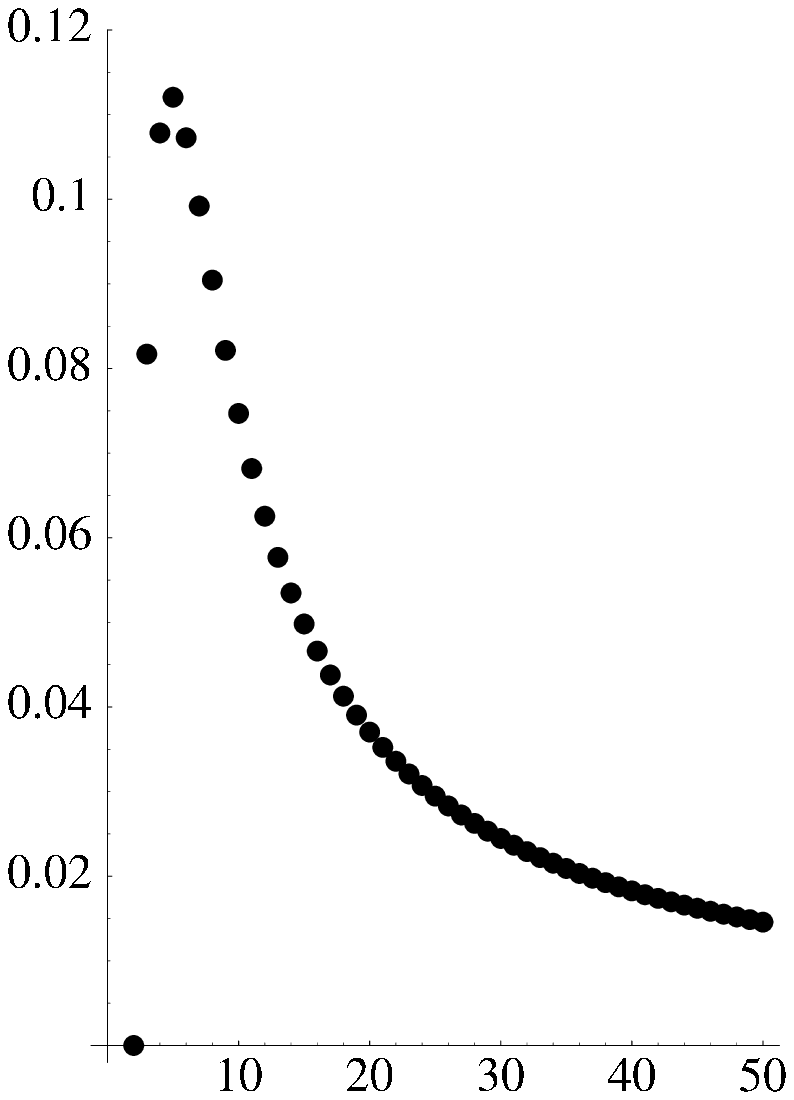}
%  \hspace{1.5cm}
%  \includegraphics[height=4cm]{PNG/h.png}
  \caption{The values $\tau_k$ for $k\ge 1$ and $1-H_r$ (see Proposition~\ref{prop:spanning}) for $r\ge 2$.}
  \label{fig:tau}
\end{figure}

\begin{proof}
  Let us first consider the case $n(d)\le 2^{(\tau_k-\varepsilon)d}$.
  We adopt the notation introduced in order to describe the first
  random model; in particular, $P_1=P=\convOp W$.  Since~$r=k+1$ is
  constant, $S$ (from the second random model) will consist of $k+1$
  affinely independent points with (very) high probability for
  large~$d$ (see~\cite{KKS95}). Thus, the first part of
  Theorem~\ref{thm:main} here implies
  \begin{equation*}
    \prob{\convOp S\text{ is a $k$-dimensional face of }P_2}\ =\ 1-\littleO{1}\ .
  \end{equation*}
  Due to~(\ref{eq:intro:neu:1}) (note $\tau_k\le\tfrac{1}{2}$), this yields 
  $$
  \prob{\convOp T\text{ is a $k$-dimensional face of }P_1}
  \ =\ 1-\littleO{1}
  $$
  for~$T$
  chosen
  uniformly at random from the $(k+1)$-subsets of (the random
  $n(d)$-set)~$W$.
  But this probability obviously is a lower bound for
  $\expect{\varphi_k(P_1)}$,
  which  proves the first part of the theorem.
  
  Now, we consider the case $n(d)\ge 2^{(\tau_k+\varepsilon)d}$.
  Similarly to the first case, the second part of
  Theorem~\ref{thm:main} here implies
  \begin{equation}
    \label{eq:neu:5}
    \probCondBig{P_2\cap \affOp S\text{ is a face of }P_2}{|S|=k+1}\ =\ \littleO{1}\ .    
  \end{equation}
  Furthermore, it is easy to see that
  \begin{multline}
    \label{eq:neu:6}
    \probCondBig{P_2\cap \affOp S\text{ is a face of }P_2}{|S\cup X|=n(d)}\\
    \le\
    \probCondBig{P_2\cap \affOp S\text{ is a face of }P_2}{|S|=k+1}
  \end{multline}
  holds. From~(\ref{eq:neu:5}) and~(\ref{eq:neu:6}) one readily deduces
  $$
  \prob{P_1\cap \affOp T\text{ is a face of }P_1}\ =\ \littleO{1}
  $$
  for~$T$ again chosen uniformly at random from the $(k+1)$-subsets
  of~$W$.  Since the number of $k$-faces of a polytope is at most the
  number of $(k+1)$-subsets of its vertex set for which the
  intersections of their affine hulls and the polytope are faces of
  the polytope, the latter probability is an upper bound for
  $\expect{\varphi_k(P_1)}$. This proves the second part of the
  theorem.
\end{proof}

\subsection{Overview of the proof of Theorem~\ref{thm:main}}

%\begin{remark}
%  \label{rem:faceChar}
%  Let~$P$ be a polytope and~$S$ a subset of its vertices.
%  \begin{enumerate}
%  \item $\conv{S}$ is a face of~$P$ \emph{if} (and only if) there is an affine
%    halfspace containing $S$ in its boundary and $V\setminus S$ in its
%    interior.
%  \item $\conv{S}$ is not contained in the boundary $\boundaryOp P$
%    of~$P$ \emph{if} the barycenter $\frac{1}{|S|}\sum_{s\in S}s$ of
%    $S$ is contained in the (relative) interior of~$P$.
%  \item If~$P\subset\R^{d}$ is a 0/1-polytope, and~$F(S)$ is the
%    smallest face of~$\cube{d}$ containing~$S$, then $\convOp S$ is a
%    face of~$P$ \emph{if and only if} $\convOp S$ is a
%    face of~$P\cap F(S)$.
%  \end{enumerate}
%\end{remark}

The structure of the proof is as follows: First, we will (in
Section~\ref{sec:reduce}) reduce the proof of Theorem~\ref{thm:main}
to a statement (Proposition~\ref{prop:spanning}) about the event
that~$S$ is not contained in a proper face of the cube, i.e., $S$ is
\emph{spanning}. (A \emph{proper} face of a polytope is any face that
is not the entire polytope, which is considered a face of itself
here.) This statement finally is proved in Section~\ref{sec:spanning}.
There we need the results of Section~\ref{sec:membership} (for
treating the cases \emph{behind} the threshold) and
Section~\ref{sec:shallow} (for the cases \emph{below} the threshold).

We will use only basic facts from polytope theory (such as in the
proof of Theorem~\ref{thm:truemain}). Consult~\cite{Zie95} in case of
doubts -- or for background information. 

Throughout the paper, $r\in\{3,4,\dots\}$ will be a constant. 

\subsection*{Acknowledgments}

I am grateful to the Mathematical Sciences Research Institute at
Berkeley for the generous support and the excellent conditions I
enjoyed during my visit in October/November 2003, when this work was
done. I thank G\"unter M. Ziegler for comments on an earlier version
of the paper.

% --------------------------------------------------
% Reduction to the spanning case
% --------------------------------------------------

\section{Reduction to the spanning case}
\label{sec:reduce}

From now on, we stick to the second model of randomness. Thus, for
some function $n\,:\,\N\rightarrow\N$, we choose the points
$S_1,\dots,S_r,X_1,\dots,X_{n(d)-r}\in\cubeVerts{d}$ \iuar, and let
$S:=\{S_1,\dots,S_r\}$, $X:=\{X_1,\dots,X_{n(d)-r}\}$, and
$P:=\conv{X\cup S}$.  Denote by $F(S)$ the smallest face of the cube
$\cube{d}$ that contains~$S$. Clearly, $P\cap F(S)$ is a face of~$P$.
Let $d(S)$ be the dimension of $F(S)$ (i.e., $d(S)$ is the number of
coordinates where not all elements of~$S$ agree).  If $F(S)=\cube{d}$
(i.e., $d(S)=d$), then we call~$S$ \emph{spanning}.

In Section~\ref{sec:spanning}, we will prove the following result
(where $\boundaryOp$ denotes the boundary operator).
\begin{proposition}
\label{prop:spanning}
  Let $r\in\{3,4,\dots\}$ and $\varepsilon>0$.
  \begin{enumerate}
  \item If $n(d)\le 2^{(1-H_r-\varepsilon)d}$ holds for all~$d$, then we have
      $$
      \probCond{\convOp S\text{ is a face of }P}{S\text{ is spanning}}\ =\ 1-\littleO{1}\ .
      $$
  \item If $n(d)\ge 2^{(1-H_r+\varepsilon)d}$ holds for all~$d$, then we have
      $$
      \probCond{\convOp S\subseteq\boundaryOp P}{S\text{ is spanning}}\ =\ \littleO{1}\ .
      $$
  \end{enumerate}
\end{proposition}

Figure~\ref{fig:tau} illustrates the threshold values $1-H_r$. 
The aim of the current section is to show that
Proposition~\ref{prop:spanning} implies Theorem~\ref{thm:main}.

\subsection{Preliminaries}

Let~$A$ be the $r\times d$ matrix whose rows are $S_1,\dots,S_r$.
Clearly, $d(S)$ equals the number of columns of~$A$ which
are neither~$\zeroVec$ (the all-zero vector) nor~$\oneVec$ (the
all-one vector). 

The random matrix~$A$ is distributed in the same way as an $r\times
d$ matrix is distributed whose columns are chosen \iuar\ from
$\{0,1\}^r$. For $t\in\{0,1\}^r$ chosen \uar, we have
$\prob{t\not\in\{\zeroVec,\oneVec\}}=1-2^{1-r}$. 

The de~Moivre-Laplace Theorem (see, e.g., \cite[Chap.~7]{Fel68}) 
yields that, for every $\delta>0$, there is a $B_{\delta}>0$ such that 
\begin{equation}
  \label{eq:reduce:1}
  \probBig{|d(S)-(1-2^{1-r})d|\le B_{\delta}\sqrt{d}}
  \ \ge \
  1-\delta
\end{equation}
holds for all large enough~$d$.

For each $\delta>0$, define
$$
J_{\delta}(d)\ :=\ \setdefBig{j\in[d]}%
                       {|j-(1-2^{1-r})d|\le B_{\delta}\sqrt{d}}\ .
$$
Thus, by~(\ref{eq:reduce:1}) we have
\begin{equation}
  \label{eq:reduce:2}
  \prob{d(S)\in J_{\delta}}\ \ge\ 1-\delta
\end{equation}
for all large enough~$d$.

Let us  denote
$$
n(S)\ :=\ \big|\setdef{i\in[n]}{X_i\in F(S)}\big|\ .
$$

\subsection{The case $\mathbf{n(d)\le 2^{(\Tilde{\tau}_r-\varepsilon)d}}$}

From elementary polytope theory one derives
\begin{equation}
  \label{eq:reduce:neu:3}
  \convOp S\text{ is a face of }P
  \ \Leftrightarrow\
  \convOp S\text{ is a face of }P\cap F(S)\ .
\end{equation}
Let $\delta>0$ be fixed and let $j_{\min}\in J_{\delta}$ such that
\begin{multline*}
  \probCond{\convOp S\text{ is a face of }P}{d(S)=j_{\min}}\\
  =\
  \min\setdefBig{\probCond{\convOp S\text{ is a face of }P}{d(S)=j}}{j\in J_{\delta}(d)}\ .
\end{multline*}
Then we have
$$
\big|d-\frac{j_{\min}}{1-2^{1-r}}\big|\ =\ \littleO{j_{\min}}\ .
$$
We therefore obtain
\begin{eqnarray}
  \label{eq:reduce:3}
  \expectCondBig{n(S)}{d(S)=j_{\min}}
     = & 2^{j_{\min}-d}(n(d)-r) \nonumber\\
   \le & 2^{j_{\min}-d+(\Tilde{\tau}_r-\varepsilon)d+\littleO{j_{\min}}} \nonumber\\
   \le & 2^{\frac{1-2^{1-r}+\Tilde{\tau}_r-1-\varepsilon}{1-2^{1-r}}j_{\min}+\littleO{j_{\min}}}\ . 
\end{eqnarray}
The fraction in the exponent equals $1-H_r-\varepsilon'$ where
$\varepsilon':=\frac{\varepsilon}{1-2^{1-r}}>0$. By Markov's inequality, we obtain 
\begin{equation}
  \label{eq:reduce:4}
  \probCond{n(S)\le 2^{(1-H_r-\varepsilon'/2)j_{\min}}}{d(S)=j_{\min}}
  \ =\
  1-\littleO{1}\ .
\end{equation}
Proposition~\ref{prop:spanning} implies
\begin{multline*}
  \probCond{\convOp S\text{ is a face of }P\cap F(S)}{d(S)=j_{\min},n(S)\le 2^{(1-H_r-\varepsilon'/2)j_{\min}}}\\
  =\ 
  1-\littleO{1}\ . 
\end{multline*}
Together with~(\ref{eq:reduce:4}), the definition of $j_{\min}$,
and~(\ref{eq:reduce:2}), this implies
$$
\prob{\convOp S\text{ is a face of }P\cap F(S)}
\ =\ 
1-\littleO{1}\ ,
$$
which, by~(\ref{eq:reduce:neu:3}), 
proves the first part of
Theorem~\ref{thm:main}.

\subsection{The case $\mathbf{n(d)\ge 2^{(\Tilde{\tau}_r+\varepsilon)d}}$}

Again, elementary polytope theory tells us
\begin{multline}
  \label{eq:reduce:neu:1}
  P\cap\affOp S\text{ is a face of }P\\
  \Rightarrow\
  P\cap\affOp S=P\cap F(S)\text{ or }\convOp S\subseteq\boundary{P\cap F(S)}\ .
\end{multline}
We omit the calculations that are necessary to prove the following lemma.

\begin{lemma}
  \label{lem:reduce}
  Let $\alpha,\beta,\gamma>0$ with $\alpha+\beta>1+\beta\gamma$,
  $\Tilde{n}(d):=\lfloor2^{\alpha d}\rfloor$, $j(d)=\beta d
  +\littleO{d}$, and let~$F$ be any $j(d)$-dimensional face
  of~$\cube{d}$. If $X_1,\dots,X_{\Tilde{n}(d)}$ are chosen \iuar\ 
  from~$\cubeVerts{d}$, then we have
  $$
  \probBig{\big|\setdef{i\in[\Tilde{n}(d)]}{X_i\in F}\big|\ge 2^{\gamma j(d)}}
  \ =\ 
  1-\littleO{1}\ .
  $$
\end{lemma}

Now we can prove the second part of Theorem~\ref{thm:main} (using
Proposition~\ref{prop:spanning}).
Let $\delta>0$ be fixed and let $j_{\max}\in J_{\delta}$ such that
\begin{multline*}
  \probCond{\convOp S\subseteq\boundary{P\cap F(S)}}{d(S)=j_{\max}}\\
  =\
  \max\setdefBig{\probCond{\convOp S\subseteq\boundary{P\cap F(S)}}{d(S)=j}}{j\in J_{\delta}(d)}\ .
\end{multline*}
  
With $\alpha:=\Tilde{\tau}_r+\varepsilon$, $\beta:=1-2^{1-r}$, and
$\gamma:=1-H_r+\varepsilon$, one easily verifies
$\alpha+\beta>1+\beta\gamma$. Since
$j_{\max}=(1-2^{1-r})d+\littleO{d}$ we thus obtain from
Lemma~\ref{lem:reduce}
\begin{equation}
  \label{eq:reduce:7}
  \probCond{n(S)\ge 2^{(1-H_r+\varepsilon)j_{\max}}}{d(S)=j_{\max}}
  \ =\
  1-\littleO{1}\ .
\end{equation}
The second part of Proposition~\ref{prop:spanning} implies
\begin{equation*}
  \probCond{\convOp S\subseteq\boundary{P\cap F(S)}}{d(S)=j_{\max},n(S)\ge 2^{(1-H_r+\varepsilon)j_{\max}}}
  \ =\ \littleO{1}\ .
\end{equation*}
Furthermore, since $\dim(\affOp S)$ is constant, we obviously have
\begin{equation*}
  \probCond{P\cap\affOp S=P\cap F(S)}{d(S)=j_{\max},n(S)\ge 2^{(1-H_r+\varepsilon)j_{\max}}}
  \ =\ \littleO{1}\ .
\end{equation*}
Together with~(\ref{eq:reduce:7}), the definition of $j_{\max}$,
and~(\ref{eq:reduce:2}), the latter two equations even hold for the
corresponding unconditioned probabilities. Thus, we have
\begin{equation*}
  \prob{\convOp S\subseteq\boundary{P\cap F(S)}
        \text{ or }
        P\cap\affOp S=P\cap F(S)}
  \ =\ \littleO{1}\ ,
\end{equation*}
which, due to~(\ref{eq:reduce:neu:1}), 
proves the second part of
Theorem~\ref{thm:main}.

% --------------------------------------------------
% Membership probabilities
% --------------------------------------------------

\section{Membership probabilities}
\label{sec:membership}

Here, we derive (from Dyer, F\"uredi, and McDiarmid's
paper~\cite{DFM92}) suitable lower bounds on $n(d)$ that, for
specified points of~$\cube{d}$, guarantee their membership in our
random 0/1-polytopes with high probability.

For any $z\in\cube{d}$, let us define
$$
p(z)\ :=\ \frac{1}{2^d}
          \min\setdefBig{|U \cap\cubeVerts{d}|}%
                        {U \subset\R^d\text{ (closed affine) halfspace}, z\in U }\ .
$$
For each $\alpha>0$, denote
$$
\cube{d}^{\alpha}\ :=\ \setdef{z\in\cube{d}}{p(z)\ge 2^{-\alpha d}}\ .
$$
For $z=(\zeta_1,\dots,\zeta_d)\in\interiorOp\cube{d}$ (the interior
of~$\cube{d}$), define
$$
H(z)\ :=\ \frac{1}{d}\sum_{j\in[d]}h(\zeta_j)\ .
$$

From Lemmas~2.1 and~4.1 of~\cite{DFM92} one can deduce the following
fact. Let us mention that in particular
the proof of Lemma~4.1 (needed for part~(2) of Lemma~\ref{lem:DFM}) is
quite hard. It is the core of Dyer, F\"uredi, and McDiarmid's
beautiful paper.

\begin{lemma}
  \label{lem:DFM}
  Let $\alpha,\varepsilon>0$.
  \begin{enumerate}
  \item If $\Tilde{n}(d)\ge 2^{(\alpha+\varepsilon)d}$ holds for all~$d$, and
    $X_1,\dots,X_{\Tilde{n}(d)}\in\cubeVerts{d}$ are chosen \iuar, then we have
    $$
    \prob{\cube{d}^{\alpha}\subseteq\convOp\{X_1,\dots,X_{\Tilde{n}(d)}\}}\ =\ 1-\littleO{1}\ .
    $$
  \item For large enough~$d$, 
    $$
    \setdef{z\in\interiorOp{\cube{d}}}{H(z)\ge 1-\alpha+\varepsilon}
    \ \subseteq\ 
    \cube{d}^{\alpha}
    $$
    holds.
  \end{enumerate}
\end{lemma}

The following straight consequence (choose
$\alpha:=1-\beta+\varepsilon/2$) of Lemma~\ref{lem:DFM} is the key
to the proof of the second part of Proposition~\ref{prop:spanning}.

\begin{corollary}
  \label{cor:DFM}
  If $\beta>0$, $\Tilde{n}(d)\ge
  2^{(1-\beta+\varepsilon)d}$  for all~$d$, and 
  $X_1,\dots,X_{\Tilde{n}(d)}\in\cubeVerts{d}$ are chosen \iuar, then
  we have
  $$
  \probBig{\setdef{z\in\interiorOp\cube{d}}%
          {H(z)\ge\beta}\subseteq\convOp\{X_1,\dots,X_{\Tilde{n}(d)}\}}
  \ =\ 1-\littleO{1}\ .
  $$
\end{corollary}

% --------------------------------------------------
% Shallow cuts of the cube
% --------------------------------------------------

\section{Shallow cuts of the cube}
\label{sec:shallow}

This section is the heart of the proof of (the first part of)
Proposition~\ref{prop:spanning}.

For $m\in\{1,2,\dots\}$, let $A(m)$ be an $r\times M$ matrix with
$M:=(2^r-2)m$ that has as its columns~$m$ copies of each vector
$v\in\{0,1\}^r\setminus\{\zeroVec,\oneVec\}$.  This choice is
motivated by the following fact (which is, however, irrelevant in this
section): If $S_1,\dots,S_r$ are chosen \iuar\ from~$\cubeVerts{M}$,
then the multiplicity~$m$ of each vector $v\in\{0,1\}^r$ among the
columns of $A(m)$ equals the expected number of appearances of~$v$ as
a column of the matrix with rows $S_1,\dots,S_r$ --- conditioned on
the event that~$S$ is spanning.

Let $s_1,\dots,s_r\in\{0,1\}^M$ be the rows of $A(m)$, and let, for
$1\le i\le r-1$, $L(i)$ be the set of indices of columns that have
precisely~$i$ ones.  We have $|L(i)|=\binom{r}{i}m$. Denote by
$\sigma(i)$ the number of ones that any of the rows has in columns
indexed by $L(i)$ (these numbers are equal for all rows).  Obviously,
we have $\sigma(i)=\frac{i}{r}\binom{r}{i}m$.

Let $b:=(\beta_1,\dots,\beta_M)$ be the barycenter of the rows
$s_1,\dots,s_r$. For each $j\in[M]$ we thus have
$\beta_j=\frac{i(j)}{r}$, if $j\in L(i(j))$. Consequently (with the
definition of $H(\cdot)$ from Section~\ref{sec:membership}),
\begin{equation}
  \label{eq:shallow:neu:0}
  H(b)\ =\ 
  \frac{1}{M}\sum_{i\in[r-1]}\binom{r}{i}\,m\,h\big(\frac{i}{r}\big)\ =\
  H_r\ .
\end{equation}

From Section~\ref{sec:membership} (see Lemma~\ref{lem:DFM}) we know
that no hyperplane in $\R^M$ that contains~$b$ can therefore cut off
significantly \emph{less} than $2^{H_r M}$ points from
$\cubeVerts{M}$, and that there are indeed hyperplanes containing~$b$
that do also not cut off significantly \emph{more} than $2^{H_r M}$
cube vertices. However, for our purposes, it will be necessary to know
that there is a hyperplane containing not only~$b$, but even the
entire set $\{s_1,\dots,s_r\}$, and nevertheless cutting off not
significantly more than $2^{H_r M}$ cube vertices.

The next result guarantees the existence of such a
hyperplane, i.e., a certain shallow cut of the cube. Its proof will
also reveal the basic reason for the appearance of the entropy function
$h(\cdot)$: It is due to the well-known fact that, for any constant $\alpha>0$, 
\begin{equation}
  \label{eq:shallow:0}
  \sum_{p\in[\alpha q]}\binom{q}{p}\ =\ 2^{h(\alpha)q+\littleO{q}}
\end{equation}
(see, e.g., \cite[Chap.~9,Ex.~42]{GKP94}).

\begin{proposition}
\label{prop:shallow}
  There are coefficients $\alpha_1,\dots,\alpha_{r-1}\in\R$, such that
  the inequality
  \begin{equation}
    \label{eq:shallow:1}
    \sum_{i\in[r-1]}\sum_{j\in L(i)}\alpha_i\xi_j
    \ \le\
    \sum_{i\in[r-1]}\alpha_i\sigma(i)
  \end{equation}
  has at most $2^{H_r M + \littleO{M}}$ 0/1-solutions
  $(\xi_1,\dots,\xi_M)\in\{0,1\}^M$. (By construction, the 0/1-points
  $s_1,\dots,s_r$ satisfy~(\ref{eq:shallow:1}) with equality.)
\end{proposition}

\begin{proof}
  Throughout the proof, we denote the components of any vectors
  $a,l,z\in\R^{r-1}$ by $\alpha_i$, $\lambda_i$, and $\zeta_i$,
  respectively. 

  For every $a\in\R^{r-1}$ and
  $l\in\N^{r-1}$, denote by 
  $\omega_a(l)$  the number of 0/1-solutions
  to~(\ref{eq:shallow:1}) with precisely~$\lambda_i$ ones in components
  indexed by $L(i)$ and define
  $$
  \omega(l)\ :=\ \prod_{i\in[r-1]}\binom{\binom{r}{i}m}{\lambda_i}\ .
  $$
  With
  $$
  L_a\ :=\ 
    \setdefBig{l\in\N^{r-1}}%
              {\sum_{i\in[r-1]}\alpha_i\lambda_i\,\le\,
               \sum_{i\in[r-1]}\alpha_i\sigma(i)}
  $$
  we thus have
  $$
  \omega_a(l)\ =\
  \begin{cases}
    \omega(l) & \text{if } l\in L_a \\
            0 & \text{otherwise}
  \end{cases}\ .
  $$
  Consequently, the number of 0/1-points
  satisfying~(\ref{eq:shallow:1}) is precisely
  \begin{equation}
    \label{eq:shallow:2}
    \sum_{l\in L_a} \omega(l)\ .
  \end{equation}
  
  If, for some~$i$, we have $\lambda_i>\binom{r}{i}m$, then clearly
  $\omega(l)=0$. Thus, the number of nonzero
  summands in~(\ref{eq:shallow:2}) is $\bigO{m^r}$.  Below, we will
  exhibit a vector $a\in\R^{r-1}$ of (constant) coefficients
  that satisfies, with $z^{\star}:=(\sigma(1),\dots,\sigma(r-1))$,
  \begin{equation}
    \label{eq:shallow:3}
    \omega(l)\ \le\ 
    \omega(z^{\star})\,2^{\littleO{M}}
  \end{equation}
  for all $l\in L_a$.  This will eventually prove the proposition,
  since we have
  \begin{eqnarray}
    \omega(z^{\star})
     =&  \prod_{i\in[r-1]}\binom{\binom{r}{i}m}{\sigma(i)} \nonumber\\
     =&  \prod_{i\in[r-1]}\binom{\binom{r}{i}m}{(i/r)\binom{r}{i}m} \nonumber\\
     =&  \prod_{i\in[r-1]}2^{h(i/r)\binom{r}{i}m + \littleO{m}} \nonumber\\
     =&  2^{\sum_{i\in[r-1]}h(i/r)\binom{r}{i}m + \littleO{m}} \nonumber\\
     =&  2^{H_rM + \littleO{M}} \nonumber
  \end{eqnarray}
  (where the third equation is due to~(\ref{eq:shallow:0}), and for
  the the last one, see~(\ref{eq:shallow:neu:0})).
  
  We now approximate the function~$\omega(\cdot)$ by Sterling's formula
  (see, e.g.,~\cite[Eq.~(9.40)]{GKP94})
  $$
  N!\ =\ \Theta\Big(\sqrt{N}\frac{N^N}{e^N}\Big)\ .
  $$
  For simplicity, we define $M_i:=\binom{r}{i}m$. Thus
  we obtain
%  $$
%  \Omega(M^{-r})
%  \prod_{i\in[r-1]}\frac{M_i^{M_i}}%
%                        {\lambda_i^{\lambda_i}(M_i-\lambda_i)^{M_i-\lambda_i}}
%  \ \le\ \omega(k)\ \le\ 
%  \bigO{M^r}
%  \prod_{i\in[r-1]}\frac{M_i^{M_i}}%
%                        {\lambda_i^{\lambda_i}(M_i-\lambda_i)^{M_i-\lambda_i}}
%  $$
  $$
  \omega(l)
  \ \le\ 
  \bigO{M^r}\prod_{i\in[r-1]}\frac{M_i^{M_i}}%
                                  {\lambda_i^{\lambda_i}(M_i-\lambda_i)^{M_i-\lambda_i}}
  $$
  (with $0^0=1$).
  Let us define the closed box
  $$
  B\ :=\ [0,M_1]\times[0,M_2]\times\dots\times[0,M_{r-1}]\ ,
  $$
  the map
  $\eta:B\rightarrow\R$ via
  $$
  \eta(z)\ :=\ 
    \prod_{i\in[r-1]}\frac{M_i^{M_i}}%
                          {\zeta_i^{\zeta_i}(M_i-\zeta_i)^{M_i-\zeta_i}}\ ,
  $$
  and the halfspace
  $$
  U_a\ :=\
    \setdefBig{z\in\R^{r-1}}%
              {\sum_{i\in[r-1]}\alpha_i\zeta_i\,\le\,
               \sum_{i\in[r-1]}\alpha_i\sigma(i)}\ .    
  $$
  We have 
  $$
  \setdef{l\in L_a}{\omega(l)>0}\ \subseteq\ B\cap U_a\ .
  $$

  By the continuity of~$\eta$ on~$B$ it hence suffices to determine
  $a\in\R^{r-1}$ such that $\eta(z^{\star})\ge\eta(z)$ holds for
  all $z\in U_a\cap\interiorOp B$.
  Note that $z^{\star}$ itself is contained in the interior
  $\interiorOp B$ of the box~$B$, where~$\eta$ is a
  differentiable function.
  
  In fact, since $\ln(\cdot)$ is monotonically increasing, we may
  equivalently investigate the function
  $\Tilde{\eta}\,:\,\interiorOp B\rightarrow\R$ defined via
  $$
  \Tilde{\eta}(z)\ :=\ \ln\eta(z)
  = \sum_{i\in[r-1]}M_i\ln M_i 
    - \sum_{i\in[r-1]}\big(\zeta_i\ln\zeta_i+(M_i-\zeta_i)\ln(M_i-\zeta_i)\big)\ ,
  $$
  and thus find a vector $a\in\R^{r-1}$ of coefficients with
  \begin{equation}
    \label{eq:shallow:neu:2}
    \Tilde{\eta}(z^{\star})\ \ge\ 
       \Tilde{\eta}(z)\qquad\text{for 
                all $z\in U_a\cap\interiorOp B$}\ .
  \end{equation}
  
  Now we choose the vector $a\in\R^{r-1}$ to be  the
  gradient of~$\Tilde{\eta}$ at~$z^{\star}$.  One easily calculates
  $$
  \alpha_i = \ln\frac{M_i-\sigma(i)}{\sigma(i)}\ .
  $$
  
  In order to prove that, with this choice, (\ref{eq:shallow:neu:2})
  holds, let $z\in U_a\cap\interiorOp B$ be
  arbitrary ($z\not= z^{\star}$).  Define $v:=z-z^{\star}$, and
  consider the function
  $\Tilde{\eta}_{z^{\star},z}:[0,1]\rightarrow\R$ defined via
  $\Tilde{\eta}_{z^{\star},z}(t):=\Tilde{\eta}(z^{\star}+tv)$.  The
  derivative of this function on $(0,1)$ is
  \begin{equation}
    \label{eq:shallow:4}
    \Tilde{\eta}_{z^{\star},z}'(t)\ =\ 
      \sum_{i\in[r-1]}v_i\ln\frac{M_i-\sigma(i)-tv_i}{\sigma(i)+tv_i}\ .
  \end{equation}
  Consider any $i\in[r-1]$, and define
  $\varrho(t):=\frac{M_i-\sigma(i)-tv_i}{\sigma(i)+tv_i}$. If $v_i\ge
  0$, then~$\varrho(t)\le\varrho(0)$, therefore, $v_i\ln\varrho(t)\le
  v_i\ln\varrho(0)=\alpha_iv_i$. If $v_i<0$, then
  $\varrho(t)>\varrho(0)$, and thus,
  $v_i\ln\varrho(t)<v_i\ln\varrho(0)=\alpha_iv_i$. Hence, in any case
  the $i$-th summand in~(\ref{eq:shallow:4}) is at most as large as
  $\alpha_iv_i$. Therefore, we obtain
  $$
  \Tilde{\eta}_{z^{\star},z}'(t)\ \le\ \sum_{i\in[r-1]}\alpha_iv_i\ .
  $$
  
  Since $z\in U_a$, we have
  $\sum_{i\in[r-1]}\alpha_iv_i\le 0$. Thus,
  $\Tilde{\eta}_{z^{\star},z}'(t)\le 0$ for all $t\in(0,1)$. Since
  $\Tilde{\eta}_{z^{\star},z}$ is continuous on $[0,1]$, we hence
  conclude $\Tilde{\eta}(z^{\star})\ge\Tilde{\eta}(z)$.
\end{proof}  

% --------------------------------------------------
% The spanning case
% --------------------------------------------------

\section{The spanning case}
\label{sec:spanning}

Using the material collected in Sections~\ref{sec:membership}
and~\ref{sec:shallow}, we will now prove
Proposition~\ref{prop:spanning} (and thus, as shown in
Section~\ref{sec:reduce}) Theorem~\ref{thm:main}.

Towards this end, let
$S_1,\dots,S_r,X_1,\dots,X_{n(d)-r}\in\cubeVerts{d}$ be chosen according
to the probability distribution induced by our usual distribution
(choosing all points \iuar) on the event that~$S:=\{S_1,\dots,S_r\}$
is spanning.  As before, define $S:=\{S_1,\dots,S_r\}$,
$X:=\{X_1,\dots,X_{n(d)-r}\}$, and $P:=\conv{S\cup X}$.

Let~$A$ be the $r\times d$ matrix with rows $S_1,\dots,S_r$. Then~$A$ is a random matrix that
has the same distribution as the $r\times d$ random matrix~$A'$
which arises from choosing each column \iuar\ from
$\{0,1\}^r\setminus\{\zeroVec,\oneVec\}$.  Therefore, if we denote the
columns of~$A$ by $t_1,\dots,t_d\in\{0,1\}^r$, then the~$t_j$ are
(independently) distributed according to the distribution
$$
\prob{t_j=t}\ =\ \frac{1}{2^r-2}\ =:\ \pi
$$
for each $t\in\{0,1\}^r\setminus\{\zeroVec,\oneVec\}$.

Define
$$
T_r\ :=\ \{0,1\}^d\setminus\{\zeroVec,\oneVec\}\ ,
$$
and denote, for every $t\in T_r$, 
$$
J(t)\ :=\ \setdef{j\in[d]}{t_j=t}\ .
$$
Let~$m\in\N$ be the largest number such that $m\le|J(t)|$ holds
for all $t\in T_r $. For each~$t$,
choose an arbitrary subset $\Tilde{J}(t)\subseteq J(t)$ with
$|\Tilde{J}(t)|=m$.

Denote by 
$$
\Delta_{\max}\ :=\ \max\setdefBig{\big||J(t)|-\pi d\big|}%
                                 {t\in T_r }
$$
the maximal deviation of any $|J(t)|$ from its expected value $\pi d$.

From the de~Moivre-Laplace Theorem (see, e.g., \cite[Chap.~7]{Fel68})
one deduces the following for each~$t\in T_r$: For every $\gamma'>0$
there is a $C'_{\gamma'}>0$ such that
\begin{equation*}
  \probBig{\big||J(t)|-\pi d\big|\le C'_{\gamma'}\sqrt{d}}\ \ge\ 1-\gamma'
\end{equation*}
holds for all large enough~$d$.  Since $|T_r|$ is a constant, one can
even derive the following stronger result from this:
For every $\gamma>0$ there is a constant $C_{\gamma}>0$ such that
\begin{equation}
\label{eq:spanning:1}
  \probBig{\Delta_{\max}\le C_{\gamma}\sqrt{d}}
  \ \ge\ 1-\gamma
\end{equation}
holds for all large enough~$d$. 

Let us define 
$$
\Tilde{D}\ :=\ 
\bigcup_{t\in T_r }\Tilde{J}(t)
$$
and $\Tilde{d}:=|\Tilde{D}|=m(2^r-2)$.
In case of $\Delta_{\max}\le C_{\gamma}\sqrt{d}$, we can deduce
\begin{equation}
  \label{eq:spanning:3}
  \Tilde{d}
  \ \ge\ 
  d-\littleO{d}\ .
\end{equation}

\subsection{The case $\mathbf{n(d)\le 2^{(1-H_r-\varepsilon)d}}$}

Let $\Tilde{S}_1,\dots,\Tilde{S}_r$ be the canonical projections of
$S_1,\dots,S_r$, respectively, to the coordinates in~$\Tilde{D}$. Then
$\Tilde{S}_1,\dots,\Tilde{S}_r$ form a matrix $A(m)$ as defined in
Section~\ref{sec:shallow}. Denote, for each $i\in[r-1]$,
$$
\Tilde{L}(i)\ :=\ \bigcup_{t\in T_r\,:\,\oneVec^T t=i }
                          \Tilde{J}(t)\ .
$$

Due to Proposition~\ref{prop:shallow},
there are coefficients $\Tilde{a}_1,\dots,\Tilde{a}_{r-1}\in\R$ such
that the inequality
\begin{equation}
  \label{eq:spanning:5}
  \sum_{i\in[r-1]}\Tilde{a}_i\sum_{j\in\Tilde{L}(i)}\xi_j
  \ \le\
  \sum_{i\in[r-1]}\Tilde{a}_i\frac{i}{r}\binom{r}{i}m\,=:\,a_0
\end{equation}
has at most $2^{H_r\Tilde{d}+\littleO{\Tilde{d}}}$ many
0/1-solutions (and $\Tilde{S}_1,\dots,\Tilde{S}_r$ satisfy the
inequality with equality).

For each $j\in[d]$ let 
$$
a_j\ :=\ 
\begin{cases}
  \Tilde{a}_i & \text{ if }j\in\Tilde{L}(i) \\
  0           & \text{ if }j\in[d]\setminus\Tilde{D}
\end{cases}\ ,
$$
i.e., $a_1,\dots,a_d$ are the coefficients of
(\ref{eq:spanning:5}) considered as an inequality for $\R^d$. 

The inequality
\begin{equation}
  \label{eq:spanning:6}
  \sum_{j\in[d]}a_j\xi_j \ \le \ a_0
\end{equation}
is satisfied with equality by $S_1,\dots,S_r$. 

Let us, for the moment, restrict our attention to the event
$\Delta_{\max}\le C_{\gamma}\sqrt{d}$. Then (\ref{eq:spanning:6}) has at
most
$$
  2^{H_r\Tilde{d}+\littleO{\Tilde{d}}}2^{d-\Tilde{d}}  
  \ =\ 
  2^{H_r d+\littleO{d}} 
$$
solutions (due to (\ref{eq:spanning:3}). 

Define the halfspace
$$
U \ :=\ \setdefBig{(\xi_1,\dots,\xi_d)\in\R^d}{\sum_{j\in[d]}a_j\xi_j \le a_0}\ ,
$$
and let~$\boundaryOp U $ be its bounding hyperplane. Thus, we have
\begin{equation}
  \label{eq:spanning:7}
  S_1,\dots,S_r\ \in\ \boundaryOp U 
  \qquad\text{and}\qquad
  \big|U \cap\cubeVerts{d}\big|\ \le\ 2^{H_r d + \littleO{d}}
  \ .
\end{equation}

Since $n(d)\le 2^{(1-H_r-\varepsilon)d}$, 
the expected number of points from~$X$ lying in~$U $ is 
$$
\frac{2^{H_r d + \littleO{d}}}{2^d}(n(d)-r)
\ \le\ 
2^{-\varepsilon d + \littleO{d}}
\ .
$$
Therefore, by Markov's inequality, 
\begin{equation}
  \label{eq:spanning:9}
  \probCond{X\cap U =\varnothing}{\Delta_{\max}\le C_{\gamma}\sqrt{d}}
  \ = \littleO{1}
\end{equation}
From~(\ref{eq:spanning:9}) and~(\ref{eq:spanning:1}) we derive
$$
\prob{\boundaryOp U \cap P=\convOp S,X\cap U=\varnothing}
\ =\
1-\littleO{1}\ ,
$$
which proves the first part
of Proposition~\ref{prop:spanning}.

\subsection{The case $\mathbf{n(d)\ge 2^{(1-H_r+\varepsilon)d}}$}

From the remarks in the introduction, we know 
\begin{equation}
  \label{eq:spanning:2}
  \prob{|S|=r}\ =\ 1-\littleO{1}\ .
\end{equation}

Let $\gamma>0$ be fixed, and assume $|S|=r$, i.e., the points
$S_1,\dots,S_r$ are pairwise disjoint.  Denote by
$b(S)=(\beta_1,\dots,\beta_d)$ the barycenter of~$S$.  For each
$t\in T_r $ and $j\in\Tilde{J}(t)$,
we have
$$
\beta_j\ =\ \frac{\oneVec^T t}{r}\ .  
$$
If $\Delta_{\max}\le C_{\gamma}\sqrt{d}$ holds, we thus have (where
the last equation is due to~(\ref{eq:spanning:3}))
\begin{eqnarray*}
  \label{eq:spanning:4}
  H(b(S)) =& \frac{1}{d}
            \big(
              \sum_{t\in T_r }
                   m h\big(\frac{\oneVec^T t}{r}\big) 
              + \littleO{d} 
            \big)\nonumber\\
         =& \frac{1}{d}
            \big(
              \sum_{i\in[r-1]}m\binom{r}{i}h(i/r) 
              + \littleO{d} 
            \big)\nonumber\\
         =& \frac{m(2^r-2)}{d}H_r + \littleO{1} \nonumber\\
         =& (1-\littleO{1})H_r + \littleO{1}\ .\nonumber\\
\end{eqnarray*}
Hence, in this case 
$$
H(b(S))\ \ge\ H_r-\tfrac{\varepsilon}{2}
$$
holds for large enough~$d$.  Since~$H$ is continuous, there is a
neighborhood~$N$ of~$b(S)$ such that $H(x)\ge H_r-\varepsilon$ holds
for all $x\in N$. Due to $n(d)\ge
2^{(1-H_r+\varepsilon)d}$, Corollary~\ref{cor:DFM} implies
$$
\probCondBig{N\subseteq\convOp X}{|S|=r,\Delta_{\max}\le C_{\gamma}\sqrt{d}}
\ \ge\ 
1-\littleO{1}\ .
$$
Together with~(\ref{eq:spanning:2}) and~(\ref{eq:spanning:1}), this
shows
$$
\prob{b(S)\in\interiorOp P}\ =\ 1-\littleO{1}\ , 
$$
which
proves the second part of Proposition~\ref{prop:spanning}.

% --------------------------------------------------
% References
% --------------------------------------------------

\end{document}